\documentclass{article}
\usepackage{amssymb,amsfonts}
\newcommand{\twomat}[4]{\left(\begin{array}{cc}#1&#2\\#3&#4\end{array}\right)}
\newcommand{\twovec}[2]{\left(\begin{array}{c}#1\\#2\end{array}\right)}
\newcommand{\twovect}[2]{\left(\begin{array}{cc}#1&#2\end{array}\right)}

\newcommand{\C}{{\mathbb{C}}}
\newcommand{\R}{{\mathbb{R}}}
\DeclareRobustCommand\openone{\leavevmode\hbox{\small1\normalsize\kern-.33em1}}
\newcommand{\id}{\mathrm{\openone}}
\newcommand{\qed}{\hfill$\square$\par\vskip24pt} 
\newcommand{\be}{\begin{equation}}
\newcommand{\ee}{\end{equation}}
\newcommand{\bea}{\begin{eqnarray}}
\newcommand{\eea}{\end{eqnarray}}
\newcommand{\beas}{\begin{eqnarray*}}
\newcommand{\eeas}{\end{eqnarray*}}

\newtheorem{theorem}{Theorem}

\newtheorem{corollary}{Corollary}

\begin{document}
\title{A Singular Value Inequality for Heinz Means}
\author{Koenraad M.R.\ Audenaert \\ Institute for Mathematical Sciences, Imperial College London, \\ 53 Prince's Gate, London SW7 2PG, United Kingdom}
\maketitle 
\textit{Abstract: }
We prove a matrix inequality for matrix monotone functions, and apply it to prove
a singular value inequality for Heinz means recently conjectured by X.\ Zhan.

\section{Introduction}
Heinz means, introduced in \cite{bhat}, are means that interpolate in a certain way between the arithmetic and geometric mean.
They are defined over $\R^+$ as
\be\label{eq:heinzs}
H_\nu(a,b) = (a^\nu b^{1-\nu}+a^{1-\nu}b^\nu)/2,
\ee
for $0\le\nu\le 1$.
One can easily show that the Heinz means are ``inbetween'' the geometric mean and the arithmetic mean:
\be
\sqrt{ab} \le H_\nu(a,b)\le (a+b)/2.
\ee
Bhatia and Davis \cite{bhatdavis} extended this to the matrix case, by showing that the inequalities remain true
for positive semidefinite (PSD) matrices, in the following sense:
\be\label{eq:agm}
|||A^{1/2} B^{1/2}||| \le |||H_\nu(A,B)|||\le |||(A+B)/2|||,
\ee
where $|||.|||$ is any unitarily invariant norm and the Heinz mean for matrices is defined identically as in (\ref{eq:heinzs}).
In fact, Bhatia and Davis proved the stronger inequalities,
involving a third, general matrix $X$,
\be\label{eq:agmx}
|||A^{1/2} X B^{1/2}||| \le |||(A^\nu X B^{1-\nu}+A^{1-\nu}XB^\nu)/2|||\le |||(AX+XB)/2|||.
\ee

X.\ Zhan \cite{zhanlama,zhan}
conjectured that the second inequality in (\ref{eq:agm}) also holds for singular values. Namely:
for $A,B\ge0$,
\be\label{eq:zhan}
\sigma_j(H_\nu(A,B)) \le \sigma_j((A+B)/2),
\ee
is conjectured to hold for all $j$. These inequalities have been proven in a few special cases.
The case $\nu=1/2$ is known as the arithmetic-geometric mean inequality for singular values, and has been proven
by Bhatia and Kittaneh \cite{bk}.
The case $\nu=1/4$ (and $\nu=3/4$) is due to Y.\ Tao \cite{tao}.
In the present paper, we prove (\ref{eq:zhan}) for all $0\le\nu\le1$.
To do so, we first prove a general matrix inequality for matrix monotone functions (Section 3).
The proof of the Conjecture is then a relatively straightforward application of this inequality (Section 4).

\bigskip

\textit{Remark:} One might be tempted to generalise the first inequality in (\ref{eq:agm})
to singular values as well:
\be
\sigma_j(A^{1/2} B^{1/2}) \le \sigma_j(H_\nu(A,B)).
\ee
These inequalities are false, however.
Consider the following PSD matrices (both are rank 2):
$$
A=\left(\begin{array}{ccc}2&4&2\\4&8&4\\2&4&4\end{array}\right),\quad
B=\left(\begin{array}{ccc}5&0&4\\0&0&0\\4&0&4\end{array}\right).
$$
Then $\sigma_2(A^{1/2} B^{1/2}) > \sigma_2(H_\nu(A,B))$
for $0<\nu<0.13$.
\section{Preliminaries}
We denote the eigenvalues and singular values of a matrix $A$ by
$\lambda_j(A)$ and $\sigma_j(A)$, respectively. We adhere to the convention that
singular values and eigenvalues (in case they are real) are sorted in non-increasing order.

We will use the positive semidefinite (PSD) ordering on Hermitian matrices throughout, denoted
$A\ge B$, which means that $A-B\ge 0$.
This ordering is preserved under arbitrary conjugations: $A\ge B$ implies $XAX^*\ge XBX^*$ for arbitrary $X$.

A matrix function $f$ is \textit{matrix monotone} iff it preserves the PSD ordering, i.e.\ $A\ge B$ implies $f(A)\ge f(B)$.
If $A\ge B$ implies $f(A)\le f(B)$, we say $f$ is \textit{inversely matrix monotone}.
A matrix function $f$ is \textit{matrix convex} iff for all $0\le \lambda\le 1$ and for all $A,B\ge 0$,
$$
f(\lambda A+(1-\lambda)B)\le \lambda f(A)+(1-\lambda)f(B).
$$

Matrix monotone functions are characterised by the integral representation \cite{bhatia,zhan}
\be
f(t) = \alpha+\beta t + \int_0^\infty \frac{\lambda t}{t+\lambda}\, d\mu(\lambda),
\ee
where $d\mu(\lambda)$ is any positive measure on the interval $\lambda\in[0,\infty)$, $\alpha$ is a real scalar
and $\beta$ is a non-negative scalar.
When applied to matrices, this gives, for $A\ge0$,
\be\label{eq:mono}
f(A) = \alpha\id+\beta A+ \int_0^\infty \lambda A(A+\lambda\id)^{-1}\, d\mu(\lambda).
\ee

The primary matrix function $x\mapsto x^p$ is
matrix convex for $1\le p\le 2$,
matrix monotone and matrix concave for $0\le p\le 1$,
and inversely matrix monotone and matrix convex for $-1\le p\le 0$ \cite{bhatia}.
\section{A matrix inequality for matrix monotone functions}
In this Section, we present the matrix inequality that we will use in the next Section
to prove Zhan's Conjecture.

\begin{theorem}
For $A,B\ge0$, and any matrix monotone function $f$:
\be\label{eq:monotone}
A f(A)+B f(B) \ge \left(\frac{A+B}{2}\right)^{1/2}\,(f(A)+f(B))\,\left(\frac{A+B}{2}\right)^{1/2}.
\ee
\end{theorem}
\textit{Proof.}
Let $A$ and $B$ be PSD.
We start by noting the matrix convexity of the function $t\mapsto t^{-1}$.
Thus
\be
\frac{A^{-1}+B^{-1}}{2}\ge\left(\frac{A+B}{2}\right)^{-1}.
\ee
Replacing $A$ by $A+\id$ and $B$ by $B+\id$,
\be\label{eq:p1}
(A+\id)^{-1} + (B+\id)^{-1} \ge 2(\id+(A+B)/2)^{-1}.
\ee
Let us now define
$$
C_k := \frac{A^k}{A+\id} + \frac{B^k}{B+\id},
$$
and
$$
M := (A+B)/2.
$$
With these notations, (\ref{eq:p1}) becomes
\be
C_0 \ge 2(\id+M)^{-1}.
\ee
This implies
\be\label{eq:p2}
C_0+\sqrt{M}C_0\sqrt{M}\ge 2(\id+M)^{-1} + 2\sqrt{M}(\id+M)^{-1}\sqrt{M} = 2\id,
\ee
where the last equality follows easily because all factors commute.

Now note: $C_k+C_{k+1} = A^k+B^k$. In particular, $C_0+C_1=2\id$, and thus (\ref{eq:p2}) becomes
\be
\sqrt{M}(2\id-C_1)\sqrt{M}\ge C_1.
\ee
Furthermore, as $C_1+C_2=2M$, this is equivalent with
\be
C_2\ge \sqrt{M} C_1 \sqrt{M},
\ee
or, written out in full:
\be
\frac{A^2}{A+\id} + \frac{B^2}{B+\id} \ge
\left(\frac{A+B}{2}\right)^{1/2}\, \left(\frac{A}{A+\id}+\frac{B}{B+\id}\right)\,\left(\frac{A+B}{2}\right)^{1/2}.
\ee

We now replace $A$ by $\lambda^{-1}A$ and $B$ by $\lambda^{-1}B$, for $\lambda$ a positive scalar.
Then, after multiplying both sides with $\lambda^2$, we obtain that
\be
\frac{\lambda A^2}{A+\lambda\id} + \frac{\lambda B^2}{B+\lambda\id} \ge \left(\frac{A+B}{2}\right)^{1/2}\,
\left(\frac{\lambda A}{A+\lambda\id}+\frac{\lambda B}{B+\lambda\id}\right)\,\left(\frac{A+B}{2}\right)^{1/2}
\ee
holds for all $\lambda\ge0$.
We can therefore integrate this inequality over $\lambda\in[0,\infty)$ using any positive measure $d\mu(\lambda)$.

Finally, by matrix convexity of the square function, $((A+B)/2)^2\le (A^2+B^2)/2$ \cite{bhatia,zhan},
we have, for $\beta\ge0$,
\be
A(\alpha\id+\beta A) + B(\alpha\id+\beta B) \ge \left(\frac{A+B}{2}\right)^{1/2}\,
\left(2\alpha\id+\beta(A+B)\right)\,\left(\frac{A+B}{2}\right)^{1/2}.
\ee
Summing this up with the integral expression just obtained, and recognising representation (\ref{eq:mono})
in both sides finally gives us (\ref{eq:monotone}).
\qed

Weyl monotonicity, together with the equality $\lambda_j(AB)=\lambda_j(BA)$, immediately yields
\begin{corollary}
For $A,B\ge0$, and any matrix monotone function $f$:
\be
\lambda_j(A f(A)+B f(B)) \ge \lambda_j\left(\frac{A+B}{2}\,(f(A)+f(B))\right).
\ee
\end{corollary}
\section{Application: Proof of (\ref{eq:zhan})}
As an application of Theorem 1 we now obtain the promised singular value
inequality (\ref{eq:zhan}) for Heinz means, as conjectured by X.\ Zhan:
\begin{theorem}
For $A,B\in M_n(\C)$, $A,B\ge0$, $j=1,\ldots,n$, and $0\le s\le1$,
\be\label{eq:zh}
\sigma_j(A^s B^{1-s} + A^{1-s} B^s) \le \sigma_j(A+B).
\ee
\end{theorem}
\textit{Proof.}
Corollary 1 applied to $f(A)=A^r$, for $0\le r\le 1$, yields
\bea
\lambda_j(A^{r+1}+B^{r+1}) &\ge& \frac{1}{2} \lambda_j\left((A+B)\,(A^r+B^r)\right) \nonumber \\
&=& \frac{1}{2} \lambda_j\left(\twovec{A^{r/2}}{B^{r/2}}\, (A+B)\,\twovect{A^{r/2}}{B^{r/2}}\right) \label{eq:th1a} \\
&=& \frac{1}{2} \lambda_j\left(\twovec{A^{1/2}}{B^{1/2}}\, (A^r+B^r)\,\twovect{A^{1/2}}{B^{1/2}}\right). \label{eq:th1b}
\eea
Tao's Theorem \cite{tao} now says that for any $2\times2$ PSD block matrix $Z:=\twomat{M}{K}{K^*}{N}\ge0$ (with $M\in M_m$ and $N\in M_n$)
the following relation
holds between the singular values of the off-diagonal block $K$ and the eigenvalues of $Z$, for $j\le m,n$:
\be
\sigma_j(K) \le \frac{1}{2} \lambda_j(Z).
\ee

The inequality (\ref{eq:th1a}) therefore yields
\bea
\lambda_j(A^{r+1}+B^{r+1}) &\ge& \sigma_j\left(A^{r/2}\, (A+B)\, B^{r/2}\right) \nonumber \\
&=& \sigma_j (A^{1+r/2} B^{r/2} + A^{r/2}B^{1+r/2}).
\label{eq:th2}
\eea
Replacing $A$ by $A^{1/(r+1)}$ and $B$ by $B^{1/(r+1)}$ then yields (\ref{eq:zh}) for $s=(1+r/2)/(1+r)$,
hence for $0\le s\le 1/4$ and $3/4\le s\le 1$.

If, instead, we start from (\ref{eq:th1b}) and proceed in an identical way as above, then we obtain (\ref{eq:zh})
for $s=(r+1/2)/(1+r)$, which covers the remaining case $1/4\le s\le 3/4$.
\qed
\section{Acknowledgments}
This work was supported by The Leverhulme Trust (grant F/07 058/U),
and is part of the QIP-IRC (www.qipirc.org) supported by EPSRC (GR/S82176/0).
The author is grateful to Prof.\ X.\ Zhan for pointing out a mistake in an earlier
draft of the manuscript.


\begin{thebibliography}{99}
\bibitem{bhatia} R.~Bhatia, \textit{Matrix Analysis}, Springer, Berlin (1997).
\bibitem{bhat} R.~Bhatia, ``Interpolating the arithmetic-geometric mean inequality
and its operator version'', Lin.\ Alg.\ Appl.\ \textbf{413}, 355--363 (2006).
\bibitem{bhatdavis} R.~Bhatia, C.~Davis, ``More matrix forms of the arithmetic-geometric mean inequality'',
SIAM J.\ Matrix Anal.\ Appl.\ \textbf{14}, 132--136 (1993).
\bibitem{bk} R.~Bhatia and F.~Kittaneh, ``Notes on matrix arithmetic-geometric mean inequalities'',
Lin.\ Alg.\ Appl.\ \textbf{308}, 203--211 (2000).
\bibitem{tao} Y. Tao, ``More results on singular value inequalities'', Lin.\ Alg.\ Appl.\ \textbf{416}, 724--729 (2006).
\bibitem{zhanlama} X.~Zhan, ``Some research problems on the Hadamard product and singular values of matrices'',
Linear and Multilinear Algebra \textbf{47}, 191--194 (2000).
\bibitem{zhan} X.~Zhan, \textit{Matrix Inequalities}, LNM1790, Springer, Berlin (2002).
\end{thebibliography}
\end{document}